\documentclass[11pt]{amsart}

\usepackage{enumerate}

\usepackage{verbatim}
\usepackage[perpage,symbol]{footmisc}
\usepackage{hyperref}
\usepackage{color}
\usepackage{graphicx}
\usepackage{soul}

\newcommand{\NN}{\mathbb{N}}

\newtheorem{definition}{Definition}[section]
\newtheorem{notation}{Notation}[section]
\newtheorem{theorem}{Theorem}[section]

\theoremstyle{remark}
\newtheorem*{remark}{Remark}

\newcommand{\RR}{\mathbb{R}}
\newcommand{\AAA}{\mathcal{A}}

\newcommand{\PRESS}{P}

\numberwithin{equation}{section}

\newcommand{\vect}{ \mathbf v} 

\DeclareMathOperator{\CDS}{CDS}
\DeclareMathOperator{\Chr}{Chr}

\DeclareMathOperator{\exon}{exon}

\DeclareMathOperator{\median}{max}

\DeclareMathOperator{\hs}{max}

\newcommand{\codons}{3-mers}

\begin{document}

\title[topological pressure for DNA sequences]{Coding Sequence Density Estimation via Topological Pressure}

\author{David Koslicki}
\address{Department of Mathematics 
\\ Oregon State University  \\ Kidder Hall 354, Corvallis, OR 97330}

\email{david.koslicki@math.oregonstate.edu}

\author{Daniel J. Thompson}
\address{Department of Mathematics \\ The Ohio State University \\
100 Math Tower, 231 West 18th Avenue, Columbus, Ohio 43210,  USA}
\email{thompson@math.osu.edu}
\thanks{D.K. was  partially supported by NSF grant DMS-$1008538$. }
\thanks{D.T. was partially supported by NSF grants DMS-$1101576$ and DMS-$1259311$}

\date{\today}

\begin{abstract}
We give a new approach to coding sequence (CDS) density estimation in genomic analysis based on the 
\emph{topological pressure}, which we develop 
from a well known concept in ergodic theory. 
Topological pressure measures the `weighted information content' of a finite word, and incorporates $64$ parameters which can be interpreted as a choice of weight  for each nucleotide triplet.
We train the parameters so that the topological pressure fits the observed coding sequence density on the human genome, and use this to give \emph{ab initio} predictions of CDS density over windows of size around $66,000$bp on the genomes of Mus Musculus, Rhesus Macaque and Drososphilia Melanogaster. While the differences between these genomes are too great to expect that training on the human genome could predict, for example, the exact locations of genes, we demonstrate that our method gives reasonable estimates for the `coarse scale' problem of predicting CDS density.  

Inspired again by ergodic theory, the weightings of the nucleotide triplets obtained from our training procedure are used to define a probability distribution on finite sequences, which 
can be used to distinguish between intron and exon sequences from the human genome of lengths between 750bp and 5,000bp. At the end of the paper, we explain the theoretical underpinning for our approach, which is the theory of Thermodynamic Formalism from the dynamical systems literature.  
Mathematica 
and MATLAB 
implementations of  our method are available at \hyperref[http://sourceforge.net/projects/topologicalpres/]{http://sourceforge.net/projects/topologicalpres/}.

\end{abstract}

\maketitle

\section{Introduction} 

\renewcommand{\thefootnote}{\fnsymbol{footnote}}
\subsection*{Overview} We present a novel approach to genomic analysis using tools from the theory of thermodynamic formalism. A number of recent influential works in mathematical biology have been based on the philosophy that the methods of statistical mechanics, and dynamical systems, can give insight into biological problems \cite{Bialek2012, Mora2010, Tkacik2006, Schneidman2006}. In this spirit, we adapt tools from thermodynamic formalism (which is a well established branch of dynamical systems, developed from ideas in statistical mechanics and information theory), to the study of bioinformatics.  The principle concept that we introduce is the \textit{topological pressure} of a finite sequence, which is adapted from a well known concept in ergodic theory. 
It is a real number which is given by counting, with weights, all distinct subwords of an exponentially shorter length that appear in the original word, and is interpreted as a weighted measure of complexity of a finite sequence\footnote{See \S\ref{tp} for a precise definition, and \S \ref{interpret} for biological interpretation}. 

The structure and organization of genomes is of central concern to the study of genome biology, and determining the distribution of coding sequences is a key component of this pursuit \cite{Berna2012, Mackiewicz2010, Salzburger2009, Kowalski2009}. Furthermore, identification of gene-rich regions in eukaryotes (especially in plants) is an ongoing field of research \cite{Ksi2013, Varshney2006, Erayman2004}. The topological pressure provides a computational tool for predicting the distribution of coding sequences and identifying such gene-rich regions. Our approach is particularly suitable for the study of novel genomes where limited training data is available. 
 This is especially useful when faced with the recent aggregation of thousands of little-studied genomes (eg. Genome 10K \cite{Genome10K}). 


 The primary goals of our analysis are:  

(1) To use the topological pressure, trained on the human genome, to give \emph{ab initio} predictions of coding sequence density on other genomes (Mus Musculus, Rhesus Macaque, Drososphilia Melanogaster). 
This establishes the key practical advantage of our approach, which is that we can predict CDS density 
using only a single moderately phylogenetically distant informant genome as training data. 


(2) To use the theory of thermodynamic formalism to turn the data encoded in the parameters used in (1) into a probability distribution which can measure the coding potential of sequences of nucleotides of lengths between $750$bp and $5000$bp. 


\subsection*{Predictions for CDS density}
The \textit{coding sequence density} (or CDS density) is the probability density function given by the bin count of coding sequences 
 in non-overlapping windows of a given size. We focus on windows of size approximately $66,000$bp for reasons we describe later. This corresponds to dividing, for example, the autosomes of the human genome into roughly $40,000$ windows. The topological pressure, which depends on 64 parameters (one for each nucleotide triplet) assigns a real number to each of these windows, and we train these parameters by maximizing the correlation with the observed CDS density on a genome. 
 
After obtaining our parameters by training on the human genome, and cross-validating our results to check we are not overfitting, we give \textit{ab initio} predictions of the CDS density of Mus Musculus, Rhesus Macaque and Drososphilia Melanogaster simply by computing the topological pressure along these genomes. We find that the correlation between topological pressure (trained on the human genome) and the observed CDS density on these genomes is $0.77, 0.73$ and $0.60$ respectively. The decrease in the correlation roughly corresponds to increasing phylogenetic distance between the human genome and the target genome. 

Our predictions of CDS density can be improved by using better training data (for example, topological pressure would estimate the CDS density of Drososphilia Melanogaster very accurately if it were trained on the genome of Drosophilia Simulans), however our results emphasize that we can still make reasonable predictions of CDS density even if we are not able to train on a close relative of the target genome. This relatively low sensitivity to organism-specific genomic traits means that although our method cannot hope to predict any finer structure of a genome (for example, the exact location of genes), our technique is advantageous for the identification of regions of high CDS density for novel genomes where refined training data is unavailable.  Our approach is also suitable for \emph{ab initio} prediction on non-mammalian genomes if a suitable model genome is chosen as training data, although we do not develop this line of research here.

\subsection*{Comparison with gene-finding techniques} In the last ten years, a number of powerful and effective 
gene-finding software packages have been developed (e.g. Augustus, Contrast, Exoniphy, Genemark HMM, FGenesh, GenSCAN, 
GeneID, N-SCAN, SNAP). While these packages were not primarily designed for estimating CDS density, this information can be inferred by taking a bin count of the predicted coding sequences. These methods, which are typically based on 
Hidden Markov Models or conditional random fields, are often very effective at gene prediction on reasonably well understood genomes, although gene sensitivity/specificity and accuracy of predicted intron-exon structure is typically much lower \cite[p.333]{Yandell2012}, \cite[fig. 1]{Flicek2007}. 

 The drawback of these gene-finding methods is that they achieve only limited success on novel 
 genomes \cite{Korf2004, Yandell2012}, as they rely on parameter files which are either partially trained on the genome under study, or use detailed data from a large number of closely related informant genomes. In particular, the training procedure requires a large number of high-quality genes and error-free assemblies, and can require data that is not yet available for new genomes \cite[p. 333]{Yandell2012}, \cite[S2.2-3]{Guigo2006}, \cite[p. 577]{Guigo2005}, \cite[p. S6.2]{Carter2006}.


We investigate the predicted CDS density given by some of these methods for comparison. We use GeneID on each of the genomes we consider, and find the predictions to be comparably accurate to the predictions yielded by our method. While the first version of GeneID was developed over ten years ago, it remains widely used, and we found that it often outperformed more recent gene-finding software for estimating CDS density. We ran GENSCAN and GenemarkHMM on all three genomes, and they were outperformed by GeneID in all three cases.

We considered a selection of the most recent gene-finding software packages (N-SCAN, Exoniphy, CONTRAST) on the genomes where suitable data was available for their implementation. CONTRAST gave the best prediction over any method considered on Drosophilia Melongaster, yielding a correlation of 0.92. This is not surprising since CONTRAST utilizes $14$ informant genomes closely related to Drosophilia Melongaster (for example, Drosophilia Simulans and Drosophilia Yakuba) to make these predictions. This amount of training data would usually be unavailable for the analysis of a novel genome. We showed that Exoniphy performed very effectively on Mus Musculus, performing as effectively as topological pressure. 

Apart from these examples, we do not give a comprehensive study of the performance of these advanced gene-finding programs for estimating CDS density, but it is our expectation that they perform as well, or better, than topological pressure when good training data is available. We emphasize that the advantage of our approach 
is the possibility of predicting CDS density in situations where insufficient data is available to effectively train the leading gene-finding software packages. 

Another advantage of our approach is its simplicity and speed: the topological pressure can predict a CDS density for a genome in a matter of seconds, while \textit{ab initio} prediction programs typically take a few hours, and evidence-based methods can take weeks \cite[p.335]{Yandell2012}.

\subsection*{A probability distribution on short segments of DNA sequences} Inspired once more by the techniques of ergodic theory, we demonstrate how our parameters determine a probability distribution on finite sequences, called an \textit{equilibrium measure}. We show that this probability distribution assigns relatively large weight to sequences which are known to be exons. This property can be used to predict the coding potential of 
DNA sequences which are orders of magnitude shorter than those on which the topological pressure is trained. 
 
The equilibrium measure is a Markov measure, so this construction can be interpreted as using the topological pressure (which makes no Markovian assumption at the training stage) to produce a Markov model suitable for identifying coding sequences. The theoretical basis for this construction is the Variational Principle from \S \ref{rigour}, which shows that the equilibrium measure maximizes a certain kind of entropy. While Markov models and entropy maximization are both familiar ideas in sequence modeling \cite{durbin}, the new ingredients here are the method for obtaining the Markov model, and the interpretation of the Markov model via topological pressure as an equilibrium measure.

The development of robust techniques that detect the coding potential of short sequences is 
an important area of research \cite{Creanza2009, Fickett1992, Gao2004, Guigo1995, Lin2008, Lin2011, Saeys2007, Washietl2011} with applications to sequence annotation as well as gene prediction.  We show that our equilibrium measure is reasonably effective in distinguishing between randomly selected introns and exons of length $750$bp in the human genome. While this approach is not as effective as the powerful comparative techniques developed in, for example, \cite{Washietl2011}, our method could be useful on novel genomes. Furthermore, this result can be interpreted as evidence that our parameters are capturing the differences in distribution of 3-mers between coding sequences and non-coding sequences.

\subsection*{Layout}The layout of the paper is as follows: In \S \ref{defs}, we develop our methodology. In \S \ref{section:results}, we present the results of our analysis of topological pressure and CDS density.  In \S \ref{section:measure from parameters}, we demonstrate how the topological pressure defines a measure on finite sequences, and show that this measure can distinguish between coding sequences and non-coding sequences. In \S \ref{rigour}, we explain the theoretical basis for our approach, and give more general definitions suitable for use in future analyses.

\section{Methodology} \label{defs}
\subsection{Topological Pressure} \label{tp}
We introduce the mathematical content of our study, and then show how it can be applied to genomic analysis. The topological pressure is a well known and well studied concept in the ergodic theory of dynamical systems.  The standard version is a quantity associated to a topological dynamical system which measures the `weighted' exponential orbit complexity of the system
\cite{PP, PT, Wa}. We introduce a finite implementation of topological  pressure which can be interpreted as a measurement of weighted information content of a finite sequence.  Topological pressure is a weighted version of topological entropy, which is a parameter free quantity introduced in \cite{Koslicki2011}. Topological entropy was shown to be effective in distinguishing between intron and exon sequences \cite{Koslicki2011}. 
For ease of exposition, we state here only a special case of the definition of topological pressure, which is the one we use for our investigation of DNA sequences, and then give a series of remarks which explain why it is defined this way. We postpone the general definition of topological pressure until \S \ref{rigour}. 

We consider finite sequences on the symbols $A, C, G, T$. We use the expressions `finite sequence' and `word' synonymously. However, `subword' has a different meaning from `subsequence': a subword is a subsequence whose entries are consecutive entries of the original sequence. 
We write a word either by using sequence notation, or juxtaposition, so the sequence $(A, G, A, T, C)$ may be written simply as $AGATC$.


We weight each word of length $3$, which we think of as a nucleotide triplet, with a positive real parameter. After choosing some order for the triplets (e.g. lexicographic order), it is convenient to record these parameters in a vector 
\begin{equation} \label{vect}
\vect =(v_{AAA}, v_{AAC},  v_{AAG}, \dots \ldots, v_{TTG}, v_{TTT})
\end{equation}
with $64$ coordinates. We are free to assume  that $\vect$ is a probability vector (we explain why in \S \ref{rigour}). We define $\Phi_{\vect}$ to be the real-valued function on the collection of words of length 3 that sends a word to its corresponding entry in $\vect$. In other words, for $a_1, a_2, a_3 \in \{A,C, G, T\}$,
\begin{equation} \label{funct}
\Phi_{\vect} (a_1a_2 a_3) :=  v_{a_1 a_2 a_3}.
\end{equation}

We can use the parameters encoded in $\vect$ to induce a weight on  a word $u=(u_1, u_2, \ldots, u_n)$ of length $n \geq 3$ by the expression
\begin{equation} \label{weights}
\mbox{``Weight assigned to } u" 
= \prod _{i=1}^{n-2} \Phi_{\vect} (u_i u_{i+1}u_{i+2}). 
\end{equation}
The topological pressure of a word with respect to $\vect$, whose formal definition follows, is given by counting the number of distinct subwords of an exponentially shorter length, with weights given by the expression \eqref{weights}.
\begin{definition} \label{Pressdef}
Let $m \geq n$ and let $w = (w_1, w_2, \ldots, w_m)$ be a finite sequence where each $w_i \in \{A,C, G, T\}$. We let $SW_n(w)$ denote the set of all subwords of length $n$ that appear in $w$, that is
\[
SW_n(w) = \{ w_{i}w_{i+1} \cdots w_{i+n-1} : i \in \{1, 2, \ldots, m-n+1\} \}.
\]

Suppose that $w$ has length $m=4^n+n-1$. 
Let $\vect$ be a probability vector of the form $\eqref{vect}$. We define the topological pressure of $w$ with respect to the parameters $\vect$, denoted $P(w, \vect)$, to be
\begin{equation} \label{pressdef}
P(w, \vect)  = \frac{1}{n} \log_4 \left( \sum_{u \in SW_n(w)} \prod _{i=1}^{n-2} v_{u_i u_{i+1}u_{i+2}} \right).
\end{equation}
\end{definition}
\begin{remark}
Since $SW_n(w)$ is defined as a set (rather than a sequence), subwords are not counted with multiplicity, so the expression inside the parentheses in \eqref{pressdef} is counting the \emph{distinct} length $n$ subwords of $w$, with weights determined by the parameters $\vect$ via the expression \eqref{weights}.  
\end{remark}

\begin{remark}
The definition above only applies to words whose length are of the form $4^n+n-1$ for some $n \in \NN$, and this is the $n$ which appears in equation \eqref{pressdef}. There are obvious ways to extend the definition of topological pressure to a word of arbitrary length (e.g. by truncating or averaging), but in this paper we need only consider words whose length are of this form. In this study, we set $n=8$, so we are looking for all distinct subwords of length 8 in a window of length $4^8+7 = 65,543$.
\end{remark}

\begin{remark}
 When all entries in $\vect$ are chosen to be equal (i.e. each entry is $\frac{1}{64}$), $P(w, \vect)$ reduces to the definition of topological entropy for finite sequences due to the first named author in \cite{Koslicki2011}. The reason we take the logarithm in base $4$ in \eqref{pressdef}, and the length of form $4^n+n-1$, rather than just $4^n$, is so that  the maximum value of the topological entropy is exactly 1, and that there exist sequences on which this maximum is attained (see discussion after Definition \ref{Pressdefgen} for details).

\end{remark}

\begin{remark}
It is possible to set up topological pressure so that instead of assigning a parameter value to each $3$-mer, we assign a parameter value to each $k$-mer for some fixed $k\geq 1$ (we give the details in \S 5). We focus on $k=3$ because of the biological importance of $3$-mers in the genetic code. Furthermore, we will see that using $4^3$ parameters neither overfits nor underfits our training data. We do not expect significant improvement to the results of this paper if we considered weightings on $k$-mers with $k>3$, and we would risk overfitting the data. 
Conversely, we checked that the case $k=3$ is a better fit for the data than $k=2$.\end{remark}



 \begin{remark}
In practical situations, we must also deal with the occurrence of non-$ACTG$ symbols (e.g. $N$). We do this by only including the subwords composed entirely of the symbols $ACTG$ 
in our computation of topological pressure. This is crucial for a genome like Rhesus Macaque where entries of $N$ appear throughout the genome. For a word $w$ with only a few occurrences of $N$, this has negligible effect on our computations. On the other hand, a word $w$ with many occurrences of $N$ has low topological pressure. This effect is consistent with our application to genomic analysis, because we want the topological pressure to predict low CDS density in regions with many occurrences of $N$. Alternatively, for very accurate genome assemblies such as the human genome, we can eliminate the vast majority of non-$ACTG$ symbols by removing the telomeres and centromeres of each chromosome. We can then restrict our attention to sequences composed entirely of $ACTG$ without difficulty.
\end{remark} 


\subsection{High topological pressure sequences: biological interpretation} \label{interpret} The sequences for which the topological pressure is large are those that balance high complexity against high frequency of $3$-mers with relatively large parameter values. This intuition is made precise by the variational principle for topological pressure from ergodic theory which we discuss in \S \ref{sec:vp}. 
Regions containing a large number of coding sequences will tend to have a different distribution of $3$-mers from those regions that do not, and we search for parameter values so that the topological pressure can detect this difference. 

It is crucial that topological pressure maximizes complexity and frequency of strongly weighted $3$-mers simultaneously: maximizing only complexity would favor random sequences, while maximizing only  the frequency of strongly weighted $3$-mers would favor sequences with very low complexity, neither of which we would expect to see in regions of high CDS density. On the other hand, we demonstrate that topological pressure, which balances both these effects, can be trained so that high topological pressure correlates with high CDS density.

Heuristically, we think of the $3$-mers which receive a relatively large parameter value in $\vect$ to be those which are sending a strong signal that we are in a coding region, while those with relatively small parameter value are those that are associated with non-coding regions, or do not send us a strong signal in either direction. 

While this heuristic may seem simplistic given the complexity of the relationship between nucleotide composition and the structure of genes, it is supported by a number of results in this paper. In $\S \ref{other app}$, we show that if we choose parameters based on this hueristic (by basing the parameters on the frequency of $3$-mers in exons), then topological pressure correlates positively with CDS density. This correlation is significantly weaker than that obtained by our training procedure, which is consistent with our expectations. Also in keeping with this heuristic, the results of $\S \ref{applofeq}$ show that the parameters obtained by our training procedure can be used to define a measure which classifies introns and exons.




\subsection{Topological pressure and CDS density estimation} \label{section:application to the human genome}
The \textit{coding sequence density} (or CDS density) is the probability density function representing the percentage of coding sequences 
 in non-overlapping windows of a given size. We describe our methodology for training the topological pressure to match the observed distribution of coding sequences on the human genome, and on other data sets.




We utilize the NCBI hg18 build 36.3 with coding sequences defined by NCBI RefSeq genes and accessed via the UCSC table browser \cite{Karolchik2004}. 
We choose a chromosome and fix an integer window size $m$ to divide the chromosome into non-overlapping windows of length $m$. The selection of the window size exhibits the typical trade-off between sensitivity and specificity: a smaller window size gives finer information on the CDS distribution, but exhibits a higher sensitivity to fluctuations in nucleotide composition. The most suitable window sizes for comparison with the topological pressure are those of the form $m=4^n+n-1$.  We focus on a window size of $65,543$ ($n=8$), as this seems to achieve a good balance. This corresponds to dividing the autosomes of the human genome into roughly $40,000$ non-overlapping windows. 
We remove any windows with non-$ACTG$ symbols, as the vast majority of these correspond to telomeres and centromeres.
 We could also carry out our analysis with different window sizes. The case $n=7$, which gives window size $m=16390$, would also be a reasonable choice and could give finer results, although it would be more computationally intensive and susceptible to noise. 
\begin{notation}
We divide each chromosome of the human genome into non-overlapping windows of length $m = 65,543$, assuming the chromosome is read in the $p$ to $q$ direction. 

Let $\Chr(i)$ denote the word which represents the $i^{th}$ chromosome of the human genome, and 
$\Chr(i, [n,m])$ denote the subword which starts at position $n$ and ends at position $m$. 
Let $w(i;n)$ denote the sequence which represents the $n^{th}$ such window along the $i^{th}$ chromosome of the human genome.\footnote{We are left with a shorter window at the end of each chromosome,  and we omit these from our study.} In other words, 
\begin{equation} \label{window}
w(i;n) = \Chr(i, [(n-1)m+1, (n-1)m+m ]).
\end{equation}
\end{notation}

\begin{definition}
We define the bin count for coding sequences in each window as follows: 
\begin{align}
\#CS(i;n):= \#  \{& { \rm RefSeq\ coding\ sequences\ with\ initial \  nucleotide} \notag \\
& {\rm  contained\ in\ } w(i;n)  \} \notag.
\end{align}
The coding sequence density on chromosome $i$ is defined to be
\[
\CDS(i,n) := \#CS(i,n)/\#CS(i),
\]
where $\# CS(i) := \# \{ {\rm Known\ coding\ sequences\ in\ }\Chr(i)\}$.
\end{definition}
For fixed $i$, $\CDS(i,n)$ is a probability density function of $n$. Note that our notation suppresses our choice of window size, as this stays fixed at $m=65,543 = 4^8+7$ throughout this work.




\begin{notation}
Given a probability vector $\vect$ with $64$ entries, as described at \eqref{vect}, we consider the topological pressure with respect to $\vect$ of each of the sequences $w(i;n)$ using the following notation:
$$\PRESS(i,n,\vect):=P(w(i;n),\vect), $$
where $P(\cdot,\cdot)$ is the topological pressure given by \eqref{pressdef}. Thus, $\PRESS(i,n,\vect)$ is the topological pressure with respect to $\vect$ of the sequence which arises as the $n^{th}$ non-overlapping window of length $65,543$ along the $i^{th}$ chromosome of the human genome.

\end{notation}
On each chromosome, i.e. for each fixed $i$, we can consider $\CDS(i,n)$ and $\PRESS(i,n,\vect)$ as functions in $n$. In fact, we want to consider these functions as $i$ ranges over a specified collection of chromosomes, most often the collection of all autosomes of the human genome. That is, the indices $i$ and $n$ are replaced with a new index $t=t(i,n)$ which tells us which window of this data set is under consideration. We modify the normalization of the coding sequence density so that $\CDS(t)$ is  a probability density function of $t$, and we consider  $\CDS(t)$ and $\PRESS(t, \vect)$ as functions in $t$. This is essentially equivalent to considering the concatenation of all the autosomes as a single sequence. Similarly, we can consider $\CDS(t)$ and $\PRESS(t, \vect)$ ranging over even larger data sets, for example  by concatenating all the autosomes from a number of different model species into a single sequence.

After fixing our data set, we train the parameters $\vect$ for maximum positive correlation between $\CDS(t)$ and $\PRESS(t, \vect)$. Our focus is mainly on the case when the data set is all autosomes of the human genome, although other data sets, both larger and smaller, are investigated where appropriate in this study. We demonstrate that our training procedure neither underfits nor overfits this training data.


\subsection{Details of training procedure} \label{training} 
For a fixed collection of chromosomes as described above, we use the Nelder-Mead \cite{Nelder1965} method to maximize the correlation between $\PRESS(t,\vect)$ and $\CDS(t)$ with respect to probability vectors $\vect$ with $64$ entries. 


Considered as functions in $t$, both $\CDS(t)$ and $\PRESS(t,\vect)$ are inherently noisy due to random fluctuations in nucleotide composition in a given chromosome as well as due to incomplete knowledge regarding coding sequences (eg. incorrectly annotated sequences). The noise in both functions is suppressed by utilizing a Gaussian filter. The radius of the Gaussian filter is chosen so that it coincides at each $t$ with the Gaussian kernel density estimation of $\CDS(t)$. 

We checked that other standard smoothing techniques (moving medians, exponential moving averages, convolution with a smoothing kernel) lead to similar results, 
and chose the Gaussian filter for our analysis due to its simplicity and speed of implementation. 

We utilize the Nelder-Mead \cite{Nelder1965} method in MATLAB \cite{matlab} to maximize the correlation between $\PRESS(t, \vect)$ and $\CDS(t)$ with respect to  $\vect$. The precision threshold for the convergence of this heuristic maximization technique was set to $10^{-6}$ and convergence was typically achieved in 10,000 steps of the algorithm.  

We focus on the case where the training data is the collection of all human autosomes. We did not include the sex chromosomes due to the well-known differences in mutation rate, selection, gene death and gene survival between the autosomes and the sex chromosomes  \cite{Wilson2009a, Wilson2009b, Kvikstad2007a, Graves2006, Makova2004}. We denote the parameters trained on all human autosomes as $\vect_{\max}$.

\section{Results } \label{results} \label{section:results} 
Using the methodology above, we present our results on CDS prediction using the topological pressure.
\subsection{Training on the human genome}
\label{section:trainingonhumans}
Our training procedure yields parameters $\vect = \vect_{\max}$ so that $\PRESS(t, \vect)$ and $\CDS(t)$ have correlation above 0.9 across all autosomes of the human genome.   
It is not at all obvious that our training procedure should work this effectively, as we are training 64 parameters to maximize correlation over approximately 40,000 data points. That our training procedure even works gives evidence that topological pressure can detect structure in the training data.




\subsection{Cross-Validation}
\label{section:Cross-Validation}
 
Since our method yields a very high correlation between $\PRESS(t, \vect)$ and $\CDS(t)$, we must check if we are overfitting the 64 parameters in $\vect$. We performed a traditional \cite{Picard1984} 7-fold cross-validation on chromosomes 1 through 21. We randomly partitioned the chromosomes into 7 equal-size samples. Of these 7, a single sample of three chromosomes was retained as a test sample. We performed the maximization procedure outlined in section \ref{training} on the remaining 6 samples and used the resulting parameters to obtain a correlation value between the topological pressure and the test sample CDS density. An average is then taken over the 7 possible choices of test sample. We repeated this procedure 50 times. The resulting mean correlation was 0.8049 with a variance of 0.0003232.
This demonstrates that the maximization procedure outlined in \S \ref{training} is not overfitting.

\subsection{Training on multiple genomes} \label{mult}
We can also train topological pressure on multiple informant genomes.
Using the methodology of \S \ref{training}, we obtained parameters by training on the data set given by concatenating all autosomes of the human, mouse (mm9) and rat (rn4) genomes. 
These are the parameters we use when we refer to `topological pressure (trained on 3 genomes)' in the following sections. 
This is intended simply to demonstrate that topological pressure can 
incorporate information from multiple genomes, and a thorough investigation of the effectiveness of this idea 
is beyond the scope of this paper. 

\subsection{CDS density estimation on the Rhesus Macaque}
\label{section:Rhesus}
We used the parameters $\vect_{\max}$ obtained from training on the human genome and showed that the correlation of the topological pressure with the coding sequence density given by RefSeq genes over all the autosomes of the Rhesus Macaque build rheMac3 was $0.726$. We repeated the experiment using the parameters trained on 3 genomes, and obtained a very slightly improved correlation of $0.738$. 
We compare this with the predictions given by GeneMarkHMM \cite{Lukashin1998}, GeneID \cite{Blanco2002}, GENSCAN \cite{Burge1997}, and N-SCAN. 


\begin{figure}
\caption{Topological pressure (trained on the human genome), CDS density predicted by GeneID, and known CDS density on chromosome 2 of rheMac3.}
\label{figure:PressureRhesus}
\includegraphics[width=5in]{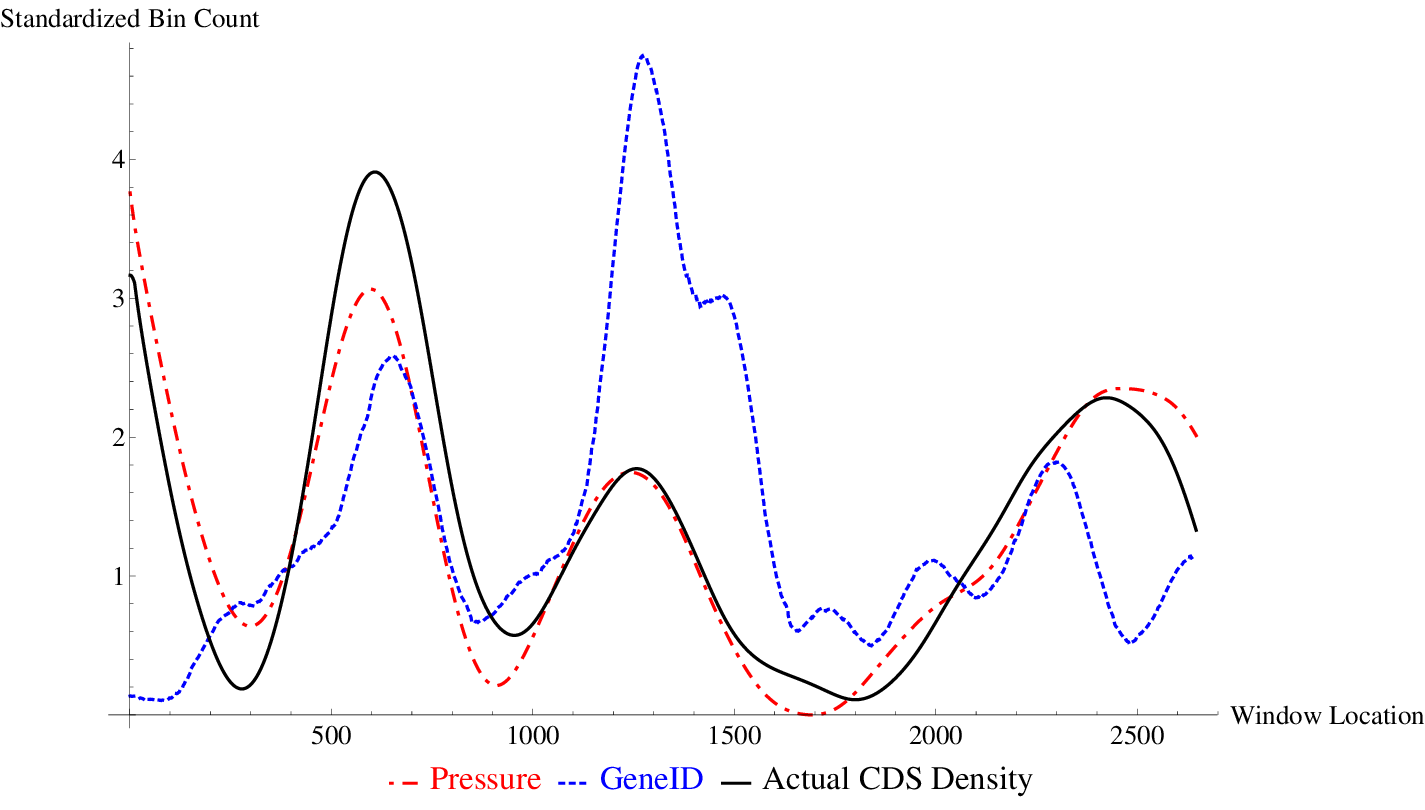}
\end{figure}


We used the GeneID and GeneMarkHMM software to obtain predicted coding sequences for the Rhesus Macaque autosomes. For GENSCAN and NSCAN, we obtained this information from the corresponding track on the UCSC table browser \cite{Karolchik2004}. For each program, we then took the bin counts of predicted coding sequences over all autosomes in the non-overlapping windows described at \eqref{window}. 
Table \ref{table:RhesusPredictionComparison} summarizes the correlation with the known coding sequence bin counts (obtained from RefSeq genes) and the bin counts predicted by each method.  Figure \ref{figure:PressureRhesus} demonstrates how well GeneID and topological pressure reconstruct the coding sequence density on chromosome 2. 

\begin{table}[!htb]
\caption{Comparison of predictions of CDS density on rheMac3.}
\begin{center}
\begin{tabular}{c|l}
\label{table:RhesusPredictionComparison}
Method & Correlation over\\ 
 & all autosomes\\ 
\hline
Topological pressure (trained on human) & 0.726\\
Topological pressure (trained on 3 genomes) & 0.738\\
GeneMarkHMM & 0.624 \\
GENSCAN & 0.402\\
GeneID & 0.660\\
N-SCAN & 0.684\\
\end{tabular}
\end{center}
\end{table}
We see that topological pressure yields the highest correlation of all the methods we looked at on this genome, and N-SCAN gave the best prediction yielded by the gene-finding programs we considered.



\subsection{CDS density estimation on Mus Musculus}
\label{section:mouse}
The correlation of the topological pressure, trained on the human genome, with the coding sequence density of the autosomes from Mus Musculus build mm9 was 0.765. We compare the topological pressure with predictions yielded by gene-finding techniques using the same methodology described in the previous section. 

We ran GeneMarkHMM on Mus Musculus genome build mm9 and obtained the GENSCAN, GeneID, and Exoniphy tracks from the UCSC table browser for this genome. 
Table \ref{table:GenePredictionComparisson} summarizes the correlation of each method with the known coding sequences density (obtained from RefSeq Genes).

\begin{table}[!htb]
\caption{Comparison of predictions of CDS density on mm9.}
\begin{center}
\begin{tabular}{c|l}
\label{table:GenePredictionComparisson}
Method & Correlation over\\
 & all autosomes\\
\hline
Top. Pressure (trained on human) & 0.765\\
GeneMarkHMM &  -0.440\\
GENSCAN & 0.695\\
GeneID & 0.817\\
Exoniphy & 0.861
\end{tabular}
\end{center}
\end{table}
Topological pressure was outperformed on this genome by GeneID and Exoniphy, but performed better than GeneMarkHMM and GENSCAN.

\subsection{CDS density estimation on Drosophila Melanogaster}
\label{section:fly}
The correlation of the topological pressure, trained on the human genome, with the coding sequence density of the autosomes from Drosophila Melanogaster build dm3 was $0.601$. This improved to $0.674$ when we used the parameters trained on $3$ genomes, and we expect that the correlation would improve significantly if we trained on a genome which was more closely related to Drosophila Melanogaster. We do not do this precisely because we want to demonstrate that we can still make reasonable predictions even when a close relative of the target genome is not available for training. 

In table \ref{table:FlyPred}, we compare the CDS prediction via topological pressure to those given by the following gene-finding techniques: GeneMarkHMM, GENSCAN, GeneID, and CONTRAST. The best performing method is CONTRAST. This may not be surprising since it uses 14 informant genomes closely related to Drosophila Melanogaster (for example, Drosophila Simulans and Drosophila Yakuba).

\begin{table}[!htb]
\caption{ Comparison of predictions of CDS density on dm3}
\begin{center}
\begin{tabular}{c|l}
\label{table:FlyPred}
Method & Correlation over\\
 & all autosomes\\
\hline
Top. Pressure (trained on human) & 0.601\\
Top. Pressure (trained on 3 genomes) & 0.674\\
GeneMarkHMM & 0.368\\
GENSCAN & 0.608\\
GeneID & 0.871\\
CONTRAST & 0.918\\
\end{tabular}
\end{center}
\end{table}

\subsection{Other approaches to parameter selection} \label{other app}
The topological pressure can be considered using parameters selected by means other than training against known data. To detect CDS density, we can select the parameters $\vect$ according to the heuristic  rule that `$3$-mers which we believe to be associated to coding sequences are assigned greater weight'. 
We give an example. 

\begin{figure}
\caption{Values of $50\times \vect_{\exon}$ and $35\times \vect_{\max}$ overlaid on the genetic code}
\label{figure:vmaxandvexon}
\vspace{5pt}
\begin{tabular}{cc}
\includegraphics[width=2.5in]{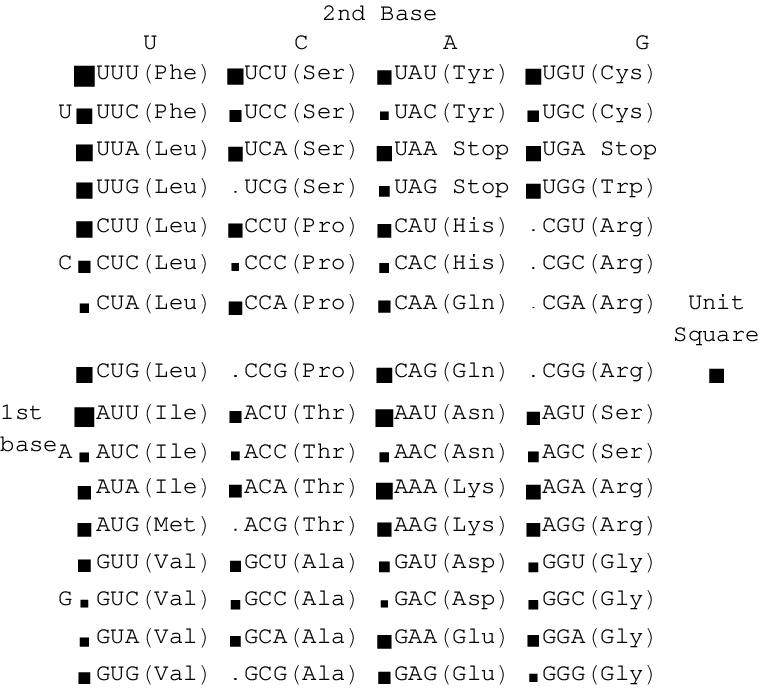} & \includegraphics[width=2.5in]{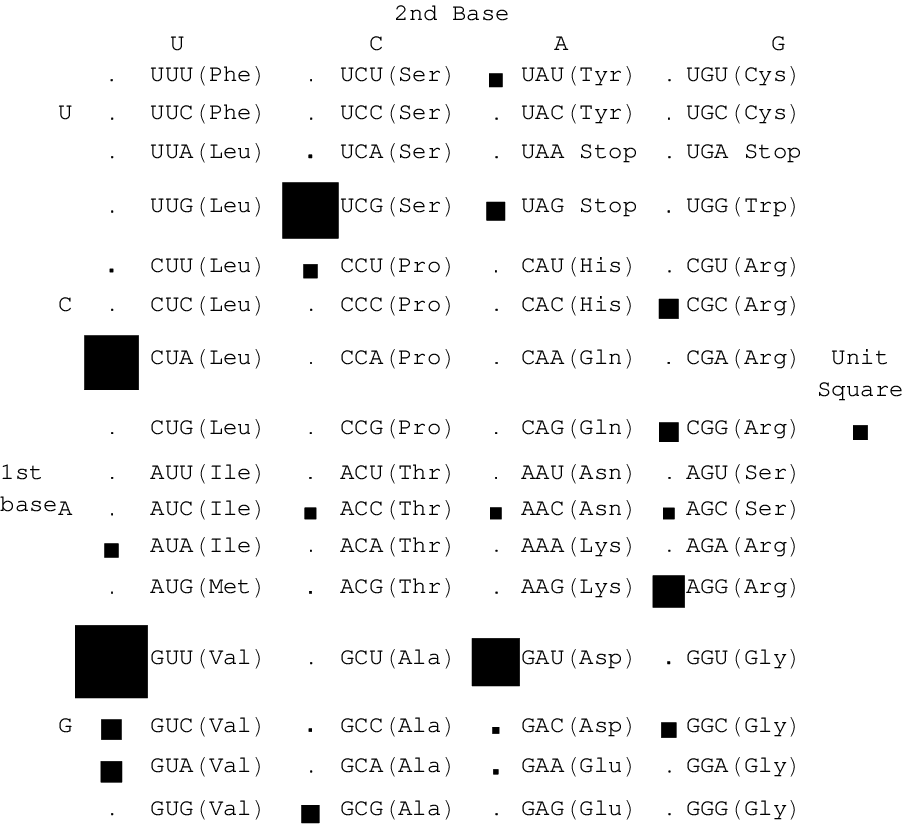}
\end{tabular}
\end{figure}

Many single sequence techniques for measuring the coding potential of DNA sequences are based upon frequencies of $n$-mers in known intronic and exonic regions \cite{Akashi2001, Comeron1998,  Creanza2009, Karlin1998}. We can use this principle to write down parameters $\vect_{\exon}$ which are based simply on the frequency of codons in the exon sequences.  More precisely, for a codon $w$, the corresponding parameter value in $\vect_{\exon}$ is assigned by the following procedure: for a segment of an autosome that corresponds to a known exon region, we count the number of times (counting overlaps) that $w$ appears, and then we sum this over all such segments. We normalize by the total number of codons (counting overlaps) that appear in the collection of segments considered, and this yields the entry in $\vect_{\exon}$ for the codon $w$. See figure \ref{figure:vmaxandvexon}.

The correlation between $\PRESS(t, \vect_{\exon})$ and $\CDS(t)$ is $0.4886$. The positive correlation matches our expectations, but it is much weaker than the correlation obtained using $\vect_{\max}$. 

\subsection{Analysis of parameter values} \label{analysis} 
\begin{figure}
\caption{Values of $35\times \vect_{\max}$ overlaid on the genetic code}
\vspace{5pt}
\label{figure:vmax}
\includegraphics[width=5in]{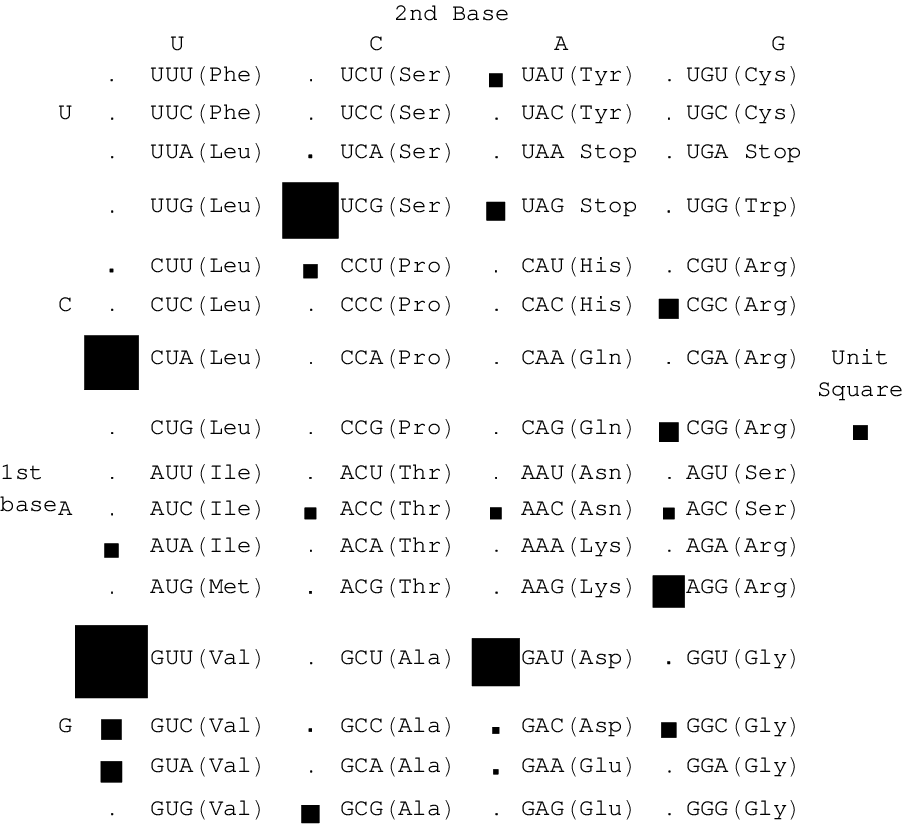}
\end{figure}

The biology enters our machinery via our choice of parameters. Since we train against known CDS density, the parameters reflect the relationship between the distribution of $3$-mers and the distribution of coding sequences along the genome. Although our method is entirely combinatorial, it would be desirable to give biological interpretation to the values assigned to $3$-mers by $\vect_{\max}$. Obvious questions include:

1) What relationship between $3$-mers and coding sequences does topological pressure really detect? We are not simply detecting the average frequency of appearance of $3$-mers in coding sequences, since the values associated to the $3$-mers by  $\vect_{\max}$ have a different, and much less uniform, distribution than average frequencies would suggest (see figure \ref{figure:vmax}). The parameters are detecting a more sophisticated relationship between the appearance of 3-mers, and their role in coding sequence formation than simply calculating frequencies, and it would be desirable to identify what biological mechanisms explain our parameter values. 

2) Do the values of $\vect_{\max}$ tell us anything about codon usage in the human genome? If we train on different genomes, what are the differences between the parameters obtained? Can this help us understand differences in codon usage between species?

A parameter sensitivity analysis will be a crucial first step in the investigation and interpretation of $\vect_{\max}$, and we hope to address these questions in future work.





We mention a feature of $\vect_{\max}$ which does match with biological intuition: 
$3$-mers made up of a single repeating nucleotide are assigned a low value by $\vect_{\max}$. Thus, the topological pressure will assign a low value to a long sequence of single repeated nucleotides. This is consistent with the presence of repetitive elements in intergenic regions of the genome.

\section{A probability measure for detection of coding potential} \label{section:measure from parameters}
An important area of research is to develop single sequence measures that effectively distinguish between short coding sequences and short non-coding sequences \cite{Creanza2009, Fickett1992, Gao2004, Guigo1995, Lin2008, Lin2011, Saeys2007, Washietl2011}. The theory of thermodynamic formalism gives us a means of selecting a Markov measure $\mu_{\vect}$, which reflects the properties of the topological pressure with respect to the parameters $\vect$. We carry out this procedure for our parameters $\vect = \vect_{\median}$ and obtain a measure that is effective for the analysis of relatively short segments of DNA sequences. We explain the theoretical underpinning for our methodology, and generalize this construction, in \S \ref{rigour}. We demonstrate that  $\mu_{\vect}$ can distinguish between coding and non-coding sequences with a reasonably high probability of success.  The advantage of using the measure  $\mu_{\vect}$ rather than the topological pressure associated to $\vect$ is that the measure is effective in analyzing relatively short DNA sequences  ($750$bp-$5000$bp).

This represents a strategy in which large scale information (parameters obtained by considering windows of $\sim 66,000$bp along the whole human genome) can be utilized to extract information at a much smaller scale (measure of a sequence of length $750$bp-$5,000$bp).

\subsection{Construction of $\mu_{\vect}$ from $\vect$}\label{section:utility of measure} We use the parameters
$$
\vect =( v_{AAA}, v_{AAC}, v_{AAG}, \dots \ldots,  v_{TTG},  v_{TTT})
$$
to define $\mu_\vect$ as a stationary Markov measure of memory $2$. In other words, our construction gives a Markov chain whose state space is the collection of all sequences of length $2$ in the DNA alphabet, and whose transition probabilities are obtained from the parameters $\vect$ by the rule \eqref{matrix20} below. The measure $\mu_v$ is then given by the standard rule for probability of a finite path of a Markov chain. See, for example, \cite{durbin} for a standard reference for these ideas in the context of biological sequence analysis. 

More precisely, 
let $\mathcal B = \{A,C,G, T \}^2$, and enumerate $\mathcal B$ by 
\[
w_1 = AA, w_2 = AC, w_3 = AG, w_4=AT, w_5=CA, \ldots, w_{16} = TT.
\]

We now use $\vect$ to define a non-negative matrix $M$ of dimension $16$ as follows. Let $M_{ij} = \vect_w$, where if $w_i = IJ$, and $w_j = JK$, then $w=IJK$. Let $M_{ij}= 0$ if the second letter in $w_i$ is not the same as the first letter in $w_j$. 
The Perron-Frobenius theorem guarantees that there is a maximal eigenvalue $\lambda >0$ and a strictly positive vector $r$ such that
\[
Mr= \lambda r.
\] 
Now define the matrix $P$ by the equation
\begin{equation} \label{matrix20}
P_{ij} = \frac{M_{ij} r_j}{\lambda r_i}.
\end{equation}
It is a standard exercise to check that $P_{ij}$ is a stochastic matrix and that there is a unique probability vector $p$ so that $pP=p$.  More explicitly, $p_i$ is given by normalizing the vector  $ l_ir_i$, where $l$ is a strictly positive left eigenvector for $M$. For $a,b,c \in \{A, C, T, G\} $, let $p(ab) = p_i$ when $ab = w_i$, and let $P(ab, bc) = P_{ij}$ when $ab= w_i$ and $bc = w_j$.



\begin{definition} \label{eq:measure definition}
We 
define a stationary probability measure $\mu_{\vect}$ on $\AAA^n$ for any fixed $n \geq 3$, by the formula 
\[
\mu_{\vect}(x_1 \cdots x_n) = p(x_1x_2) P(x_1x_2, x_2x_3) P(x_2x_3, x_3x_4) \ldots P(x_{n-2} x_{n-1}, x_{n-1} x_n)
\]
for each $x_1 \cdots x_n \in \AAA^n$.
\end{definition}
As an illustrative example, to compute $\mu_\vect$ for the word $GCTAC$, we  use the formula
\[
\mu_{\vect}(GTCAC)=p(GT)P(GT,TC)P(TC,CA)P(CA,AC),
\]
and read off the appropriate values for the right hand side of the equation.
\subsection{Detection of coding potential using $\mu_{\vect}$} \label{applofeq}
We take the measure $\mu_{\vect}$ corresponding to the parameters $\vect = \vect_{\median}$ from \S \ref{section:results}. The construction of the measure is designed so that $\mu_{\vect}$ reflects the properties of the topological pressure with respect to $\vect$ (see \S \ref{Gibbs} for details). Thus, we expect that the sequences with relatively large measure are those with higher coding potential.  

We demonstrate this phenomena by showing that $\mu_{\vect}$ can partially distinguish between a randomly selected assortment of intron and exon sequences of length $750$bp. { Sequences of this length are produced by some next-generation sequencing platforms (e.g. PacBio RS II, Roche GS FLX+)}. We randomly select $5,000$ intron sequences and $5,000$ exon sequences from human chromosome 1, {and truncate to a length of} {750bp}. 
These sequences are completely un-preprocessed: no information such as ORF's, stop/start codons or repeat masking is utilized. 

\begin{figure}
\caption{Histogram of $\log(\mu_{\vect})$ on 5,000 Introns and Exons of length 750bp}
\label{figure:muhistogram}
\includegraphics[width=5in]{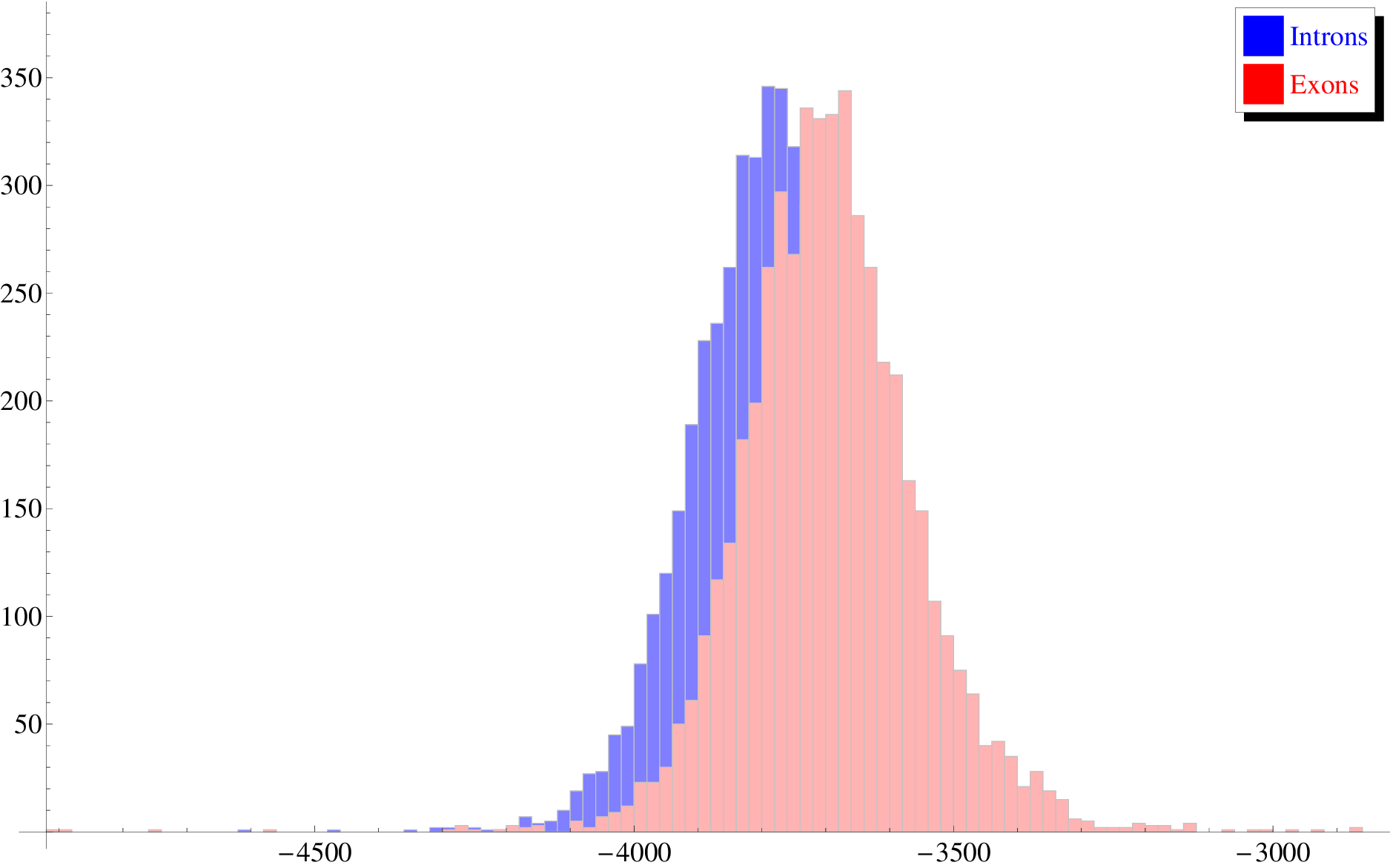}
\end{figure}

As expected, $\mu_{\vect}$ typically weights exon sequences more heavily than intron sequences. This is demonstrated by figure \ref{figure:muhistogram}, which shows the histogram of $\log(\mu_{\vect})$ evaluated on the test sequences. 
The area under the ROC (true positive rate vs. false positive rate) curve is {0.701}.

\begin{figure}
\caption{ROC curve for $\log(\mu_{\vect})$ on 5,000 Introns and Exons of length 5,000bp}
\vspace{5pt}
\label{figure:roc}
\includegraphics[width=3in]{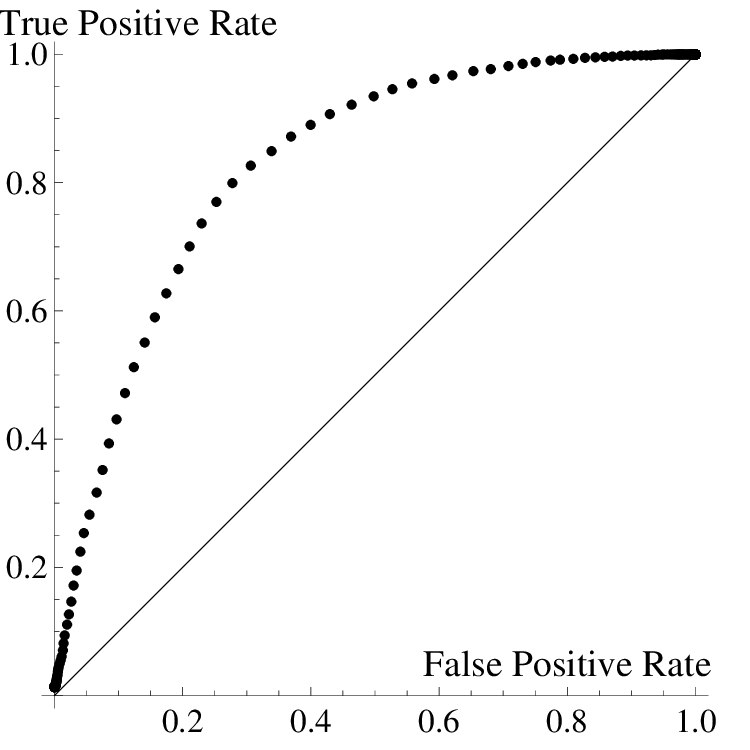}
\end{figure}

We repeated the experiment for a randomly selected assortment of introns and exons of length {$5,000$bp}, and include in figure \ref{figure:roc} the ROC curve associated to the resulting $\mu_{\hs}$. The area under the ROC curve {increased to 0.826}. 

We expect that this classification could be improved, particularly for shorter sequences, by the following strategies:

1) considering more parameters in the topological pressure, which would yield a Markov measure of higher order (as described in \S \ref{section:potentialtomeasure});

2) training on windows of much smaller length than the $\sim 66,000$bp used previously.

We do not pursue this here, and as it stands, the comparative techniques already available on the human genome \cite{Washietl2011, Creanza2009} are more accurate classifiers of introns and exons than $\mu_{\vect}$.  
Nevertheless, the equilibrium measure could potentially be a useful classifier of introns and exons on less well understood genomes. Furthermore, these results demonstrate how the parameter values for topological pressure can be used to construct a Markovian model, which captures the biological information incorporated into our machinery via the training data.

\section{Theoretical Underpinnings} \label{rigour}
Topological pressure and equilibrium measures are the principle object of study of thermodynamic formalism, which is a well established branch of ergodic theory and dynamical systems. Standard references are \cite{Ba, Bo, PP, PT, Wa}. In this section, we explain the connections between the present work and the classical theory. 

First, we extend the definition of topological pressure for finite sequences to full generality. 
Let  $\mathcal A$ be an alphabet, that is, a finite collection of symbols, 
and $| \mathcal A |$ denote the number of elements in $\mathcal A$. We denote the space of sequences of length $n$ by $\AAA^n$,   the space of finite sequences (of any length) $\AAA^{<\NN}$, the space of finite sequences of length at least $n$ by  $\mathcal A^{\geq n}$ 
and the space of infinite sequences by   $\Sigma = \mathcal A^\NN$.  For a suitable choice 
of $k$, we select a weight for each word in $\AAA^k$. The weights can be encoded by a vector $\vect$, as in \S \ref{defs}, or by a function $\psi: \AAA^k \mapsto \RR$ so that $\psi(w)$ is the weight assigned to $w$. We use the latter notation here, because it is consistent with the conventions of the dynamical systems literature. In ergodic theory, $\psi$ is customarily called the `potential function'. We avoid this terminology as the word `potential' has other meanings in biology. Often, for a function $\psi >0$, we are interested in the weights corresponding to $\Phi = \log \psi$. 

For $n \geq k$, we assign a weight to each word $u \in \AAA^{n}$ by the rule 
\begin{align*}
\mbox{``Weight assigned to } u" &= \exp \left \{ \sum_{i=1}^{n-k+1} \psi (u_iu_{i+1} \cdots u_{i+k-1}) \right \}\\ 
&= \prod _{i=1}^{n-k+1} \Phi(u_iu_{i+1} \cdots u_{i+k-1}), \mbox{ if } \Phi>0, \psi = \log \Phi.
\end{align*}
\begin{definition} \label{Pressdefgen}
Let $\psi: \AAA^k \mapsto \RR$, $m \geq n \geq k$ and let $w = (w_1, w_2, \ldots, w_m)$ be a finite sequence where each $w_i \in \AAA$. We let $SW_n(w)$ denote the set of all subwords of length $n$ that appear in $w$, that is
\[
SW_n(w) = \{ w_{i}w_{i+1} \cdots w_{i+n-1} : i \in \{1, 2, \ldots m-n+1\} \}.
\]
Now suppose that $w$ has length $4^n+n-1$, i.e. suppose $m = 4^n+n-1$. Then we can define the topological pressure of $w$ with respect to $\psi$, denoted $P(w, \psi)$, to be
\begin{equation} \label{pressdefgen}
P(w, \psi)  = \frac{1}{n} \log_{ | \AAA |}\left( \sum_{u \in SW_n(w)} \exp \left \{ \sum_{i=1}^{n-k+1} \psi (u_iu_{i+1} \cdots u_{i+k-1}) \right \} \right).
\end{equation}
If $\Phi >0$, and $\psi = \log \Phi$, where $\log$ denotes natural logarithm, then
\begin{equation} \label{pressdefgen}
P(w, \log \Phi)  = \frac{1}{n} \log_{ | \AAA |}\left( \sum_{u \in SW_n(w)} \prod _{i=1}^{n-k+1} \Phi(u_iu_{i+1} \cdots u_{i+k-1}) \right).
\end{equation}
For a word $w$ with $|\AAA|^n+n-1 \leq |w|< \AAA^{n+1} +n$, we define the topological pressure of $\psi$ on $w$ to be 
the topological pressure of $\psi$ on the first $|\AAA|^n+n-1$ symbols of $w$.

\end{definition}
Definition \ref{Pressdefgen} generalizes Definition \ref{Pressdef} because $P(w, \vect) = P(w, \log \Phi_{\vect})$,  where $\Phi_{\vect}$ is the function defined at \eqref{funct}. When $\psi =0$, \eqref{pressdefgen} reduces to the definition of topological entropy for finite sequences due to the first named author in \cite{Koslicki2011}. We denote the greatest topological pressure for words of length $4^n+n-1$ by
\begin{equation}
P_{\max} (n, \psi) = \max \{ P(w, \psi) : |w| = 4^n+n-1 \}.
\end{equation}
 For each $n$, there exists a word $w_{\max}^n$ of length $ 4^n+n-1$ which has every word of length $n$ as a subword. This follows easily from the fact that the De Brujn graph is a Hamiltonian graph, see \cite[obs. 1.6]{Gheorghiciuc2008}. It follows that  $P_{\max} (\psi, n) = P(w_{\max}^n, \psi)$, and thus $P_{\max}(n, 0) = 1$. 


Taking a multiple of $\Phi$ (equivalently adding a constant to $\psi$) does not affect the quantities associated to the topological pressure that we study in this paper, particularly correlation with the CDS density developed in \S \ref{section:application to the human genome}.
For any $t>0$, and word $w$ of length $ 4^n+n-1$ we have the formula
\begin{equation}
\label{eq:potential normalize}
P(w, \log t \Phi) = \frac{n-k}{n} \log_{|\AAA |} t + P(w, \log \Phi).
\end{equation}
Since the difference between $P(w, \log t \Phi)$ and $P(w, \log \Phi)$ is a constant independent of $w$, the correlations studied in \S \ref{defs} will remain unchanged when normalizing $\Phi$.  Hence we are free to assume that $\vect$ is a probability vector in \S \ref{defs}. 

\subsection{Equilibrium measures} \label{section:potentialtomeasure}
Given a function $\psi: \AAA^k \mapsto \RR$, there is a unique probability measure $\mu_\psi$, called the \textit{equilibrium measure} for $\psi$, whose properties reflect those of the topological pressure with respect to $\psi$.
The measure $\mu_{\vect}$ constructed in \S \ref{section:utility of measure} is an equilibrium measure. In this section, we describe how to construct equilibrium measures and explain the theoretical basis for their useful properties. 

The construction is a generalization of the construction of $\mu_{\vect}$, and a special case of more general expositions given in \cite{Ba, Bo, PP, PT, Wa}. We take our finite alphabet $\AAA$, and a function $\psi: \AAA^k \mapsto \RR$. 

  Let $\mathcal B = \AAA^{k-1}$ and enumerate $\mathcal B$ by some natural ordering.
Define a $1-0$ square matrix $S$ of dimension $|\AAA|^{k-1}$ as follows. Let $S_{ij} =1$ if and only if the word obtained by omitting the first symbol of $w_i$ is the same as the word obtained by omitting the last symbol in $w_j$. In this case, define $\pi(w_i, w_j) \in \AAA^k$ as the word $w_ib$, where $b \in \AAA$ is the last symbol in $w_j$. Equivalently, $\pi(w_i, w_j) = aw_j$, where $a \in \AAA$ is the first symbol of $w_i$. 

We now use $\psi$ to define a non-negative matrix $M$ of dimension $|\mathcal B|^2$ as follows. If $S_{ij} =1$, then let
\begin{equation} \label{matrixM}
M_{ij} = e^{\psi(\pi(w_i, w_j))},
\end{equation}
and if $S_{ij} =0$, then let $M_{ij} =0$.
The Perron-Frobenius theorem gives a maximal eigenvalue $\lambda >0$ and a strictly positive vector $r$ such that
\[
Mr= \lambda r. 
\] 
Now define a matrix $P$ of dimension $|\mathcal B|^2$  by
\begin{equation} \label{matrix2}
P_{ij} = \frac{M_{ij} r_j}{\lambda r_i}.
\end{equation}
It is easy to check that $P_{ij}$ is a stochastic matrix and that there is a unique probability vector $p$ so that $pP=p$.  More explicitly, $p_i$ is given by normalizing the vector  $ l_ir_i$, where $l$ is a strictly positive left eigenvector for $M$. 



To define a measure on $\AAA^\NN$, it suffices to define the measure on the cylinder sets 
\begin{equation}\label{cylinder}
[x_1 \cdots x_n] := \{ y \in \AAA^\NN \mid y_1 = x_1, y_2 = x_2, \ldots, y_n = x_n \}, 
\end{equation}
since these are open sets which generate the natural topology on $\AAA^\NN$ (see \cite{Wa}).

\begin{definition} \label{eq:measure definition general}
We 
define a probability measure $\mu_\psi$ on $\AAA^\NN$ 
by the formula 
\begin{equation} \label{eqmeas}
\mu_\psi([x_1 \cdots x_n]) = p_{i_1} P_{i_1i_2} P_{i_2i_3} \cdots P_{i_{n-k}i_{n-k+1}},
\end{equation}
for any $x_1 \cdots x_n \in \AAA^n$ with $n \geq k$, where $w_{i_1} = x_1 \cdots x_k$, $w_{i_2} = x_2 \cdots x_{k+1}$, $\ldots$, $w_{i_{n-k+1}} = x_{n-k+1} \cdots x_n$. We call the measure $\mu_\psi$ the equilibrium measure for $\psi$ on $\AAA^\NN$.
\end{definition}
For any fixed $n \geq k$, we can take the value assigned to each $x_1 \cdots x_n$ by the formula \eqref{eqmeas} to define a probability measure on $\AAA^n$
, which we refer to as the equilibrium measure for $\psi$ on $\AAA^n$. Thus, the probability measure $\mu_\vect$ from Definition \ref{eq:measure definition} is the equilibrium measure for $\log \Phi_\vect$ on $\{A, C, G, T\}^n$.


\subsection{Relation to theory of dynamical systems: the full shift and the Variational Principle} \label{sec:vp}
In the next few sections, we recall the classical  theory from dynamical systems which explains the importance of $\mu_\psi
$. We demonstrate the relationship between the concepts introduced in this paper and the dynamics of the full shift (defined below).
\begin{definition}
The \emph{full shift} over an alphabet $\mathcal A$ is the dynamical system $(\Sigma, \sigma)$, where $\Sigma = \mathcal A^{\NN}$ is the space of infinite sequences on $\mathcal A$, and $\sigma$ is the \emph{shift map} $\sigma: \Sigma \to \Sigma$, which is the map defined by `shifting' a sequence one position to the left. That is, for $(x_1, x_2, x_3, \ldots) \in \Sigma$,
\[
\sigma((x_1, x_2, x_3, \ldots)) := (x_2, x_3, x_4, \ldots).
\]
\end{definition}
\begin{definition}
Given a continuous function $\psi: \Sigma \to \RR$, the \emph{topological pressure}  of $\psi$ on $\Sigma$ is defined to be:
\[
P(\Sigma, \psi) = \lim_{n \to \infty} \frac 1n \log \left ( \sum_{u \in \mathcal A^n} \exp \sum_{i=0}^{n-1} \psi (\sigma^iu) \right ) .
\]
\end{definition}
The following result \cite{PT, Wa} gives the fundamental relationship between the topological pressure and $\sigma$-invariant probability measures\footnote{that is, probability measures which satisfy $\mu(\sigma^{-1} A) = \mu (A)$ for all Borel sets $A \subset \Sigma$.} on $\Sigma$.  
\begin{theorem}[Variational Principle] \label{thm:vp}
The topological pressure of $\psi$ on $\Sigma$ satisfies:
\begin{equation} \label{vp}
P(\Sigma, \psi) = \sup_{m} \left \{h_m + \int \psi d m \right \},
\end{equation}
where the supremum is taken over all $\sigma$-invariant probability measures on $\Sigma$, and 
$h_m$ denotes the measure theoretic entropy, given by
\[
h_{m} = \lim_{n \to \infty} -\frac{1}{n} \sum_{w\in \AAA_n} m([w])\log m([w]).
\]  
A measure achieving the supremum in the \eqref{vp} is called an \emph{equilibrium measure} for $\psi$. 
\end{theorem}
The following result, proved in \cite[\S 4]{PT}, tells us that the measure constructed in the previous section is indeed an equilibrium measure in this sense. 
\begin{theorem}
The measure $\mu = \mu_\psi$ defined in Definition \ref{eq:measure definition general} is the unique equilibrium measure for $\psi$ (in the sense of Theorem \ref{thm:vp}), and
\[
P(\Sigma, \psi) =  h_\mu + \int \psi d \mu = \log \lambda,
\]
where $\lambda$ is the Perron-Frobenius eigenvalue of the matrix \eqref{matrixM}.
\end{theorem}
The Variational Principle illustrates the trade-off between structure and complexity which is detected by the topological pressure, simultaneously maximizing entropy (which is itself maximized by the uniform measure) and the integral of $\psi$ (which is itself maximized by a Dirac measure). 

\subsection{The Gibbs property} \label{Gibbs}
The relationship between $\psi$ and $\mu_\psi$ is captured by the \textit{Gibbs property}, established in \cite{Bo, PP}. To simplify notation, we return to the case of $\psi: \AAA^3 \mapsto \RR$, which is the important case for this paper.
\begin{theorem}[Gibbs property] 
For $\psi: \AAA^3 \mapsto \RR$ and any $w \in \AAA^n$,
\[
\mu_{\psi}([w]) \asymp \exp \{-nP(\Sigma, \psi) + \sum_{i=1}^{n-2}\psi(w_iw_{i+1} w_{i+2}) \},
\]
where $[w]$ is the cylinder set defined at \eqref{cylinder}, and $a_n \asymp b_n$ means there exists a constant $C >1$ so that $C^{-1} \leq a_n/b_n \leq C$ for all $n$. 
\end{theorem}
Thus, if $\psi = \log \Phi$ and we normalize $\psi$ so that $P(\Sigma, \psi) = 0$ (which is done by taking a suitable multiple of $\Phi$), then
\begin{equation} \label{Gibbs2}
\mu_{\psi}([w]) \asymp \prod _{i=1}^{n-2}\Phi(w_iw_{i+1} w_{i+2}).
\end{equation}
In the context of \S \ref{applofeq}, this formula provides the intuition that sequences which have a relatively high frequency of words $w \in \AAA^3 $ where $v_w$ is large, and a relatively small frequency of words $w \in \AAA^3$  where $v_w$ is small, will be assigned relatively large measure by $\mu_\vect$.  This gives a theoretical underpinning for using $\mu_\vect$ to predict the coding potential of short sequences.

\subsection{Relationship between topological pressure for finite sequences and  topological pressure on the full shift}
We continue to focus on the case when $\psi: \{A,C, T, G\}^3 \to \RR$ for simplicity, and we write $\Sigma$ for the full shift on $\{A, C, T, G\}$. The following result is essentially that of \cite[Theorem 7.30]{Wa}. Let $M$ be the matrix constructed in \S \ref{section:utility of measure}, and recall that $\lambda$ is its Perron-Frobenius eigenvalue. We consider the matrix norm of $M$ given by $\|M\| = \sum_{i,j} |m_{ij}|$.
\begin{theorem}
We have $$P_{max}(\psi, n) = \frac{1}{n}\log_4 \left ( \sum_{u \in \mathcal A^n} \exp  \left \{\sum_{i=1}^{n-2} \psi (u_i u_{i+1} u_{i+2})  \right \} \right )  = \log_4 \|M^{n-2}\|^{1/n},$$
The sequence $\|M^{n-2}\|^{1/n}$ converges to $\lambda$ exponentially fast as $n \to \infty$.  
\end{theorem}
This theorem tells us that for large $n$, $P_{max}(\psi, n)$ is very close to $\log_4 \lambda$. Since $P(\Sigma, \psi) = \log \lambda$, this describes the relationship between topological pressure for finite sequences and  topological pressure on the full shift. 

\section{Conclusion} 
We demonstrated that the topological pressure can train 
on the human genome to fit the observed bin count of coding sequences on windows of size approximately $66,000$bp. 
We showed that topological pressure, trained on the human genome, gave effective estimates of CDS density on Rhesus Macaque, Mus Musculus and Drosophilia Melanogaster, despite the phylogenetic distance between these target genomes and the informant genome. We compared these results with predictions of CDS density yielded by a selection of current gene-finding packages. These often performed extremely well, but required detailed organism-specific training data that is not required to train the topological pressure, and is not typically available for novel genomes.


We showed that the topological pressure defines a probability measure which can distinguish between segments of human intron and exon sequences of length between 750bp and 5000bp. Finally, we established the theoretical basis for our results, adapting ideas and results from ergodic theory. 

\section*{Acknowledgements}
Portions of this work were completed while D.K. was a postdoctoral fellow at the Mathematical Biosciences Institute of the Ohio State University, and while D.K. and D.T. were members of the Mathematics Department at the Pennsylvania State University. A preliminary version of this work is included in D.K.'s PhD thesis at Penn State. The authors wish to thank the anonymous referees of this paper, whose input has greatly benefited this work.

\bibliographystyle{abbrv}
\bibliography{PressureDNA2}

\end{document}